\documentclass[12pt]{article}

\usepackage{latexsym}
\usepackage{amsmath}
\usepackage{amsfonts}
\usepackage{amssymb}

\newcommand{\bean}{\begin{eqnarray}}
\newcommand{\eean}{\end{eqnarray}}
\newcommand{\bea}{\begin{eqnarray*}}
\newcommand{\eea}{\end{eqnarray*}}
\newcommand{\bsa}{\begin{subarray}{c}}
\newcommand{\esa}{\end{subarray}}
\newcommand{\bi}{\begin{itemize}}
\newcommand{\ei}{\end{itemize}}

\newtheorem{lemma}{Lemma}[section]
\newtheorem{thm}[lemma]{Theorem}

\newtheorem{cor}[lemma]{Corollary}
\newtheorem{defn}[lemma]{Definition}
\newtheorem{propn}[lemma]{Proposition}

\def\ssp{\def\baselinestretch{1.0}\large\normalsize}

\title{ \bf Logarithmic Vector-Valued Modular Forms }
\author{Marvin Knopp \\
Geoffrey Mason\thanks{Supported by NSA and NSF}}
\date{}
\begin{document}
\ssp
\maketitle
 \begin{abstract}
 \noindent
 We consider logarithmic vector- and matrix-valued modular forms of integral weight $k$ associated with a $p$-dimensional representation $\rho: SL_2(\mathbb{Z}) \rightarrow GL_p(\mathbb{C})$ of the modular group, subject only to the condition that $\rho(T)$ has eigenvalues of absolute value $1$. The main result is the construction of meromorphic matrix-valued Poincar\'{e} series associated to $\rho$  for all large enough weights. The component functions are logarithmic $q$-series, i.e., finite sums of products of $q$-series and powers of $\log q$. We derive several consequences, in particular we
 show that the space $\mathcal{H}(\rho)=\oplus_k  \mathcal{H}(k, \rho)$ of all holomorphic logarithmic vector-valued modular forms associated to $\rho$ is a free module of rank $p$ over the ring of classical holomorphic modular forms on $SL_2(\mathbb{Z})$.

\end{abstract}

 \section{\Large \bf Introduction}
The present work is a natural sequel to our earlier articles
on `normal' vector-valued modular forms \cite{KM1}, \cite{KM2}. The component functions of a normal vector-valued modular form $F$ are $q$-series with at worst real exponents. Equivalently,
the finite-dimensional representation $\rho$ associated with $F$ has the property that
$\rho(T)$ is (similar to) a matrix that is unitary and diagonal. Here, 
$T = \left(\begin{array}{cc}1 & 1 \\0 & 1\end{array}\right)$.

\bigskip
In the case of a general representation, $\rho(T)$ is not necessarily diagonal but may always be assumed to be in Jordan canonical form\footnote{We actually  use a modified Jordan canonical form. See Subsection 2.2 for details.}. This circumstance leads to \emph{logarithmic}, or \emph{polynomial} $q$-expansions for the component functions of a vector-valued modular form associated to $\rho$ (see Subsection 2.2), which take the form
\begin{eqnarray}\label{logform1}
f(\tau) = \sum_{j=0}^t (\log q)^j h_j(\tau),
\end{eqnarray}
where the $h_j(\tau)$ are ordinary $q$-series.
There follow naturally the definition of logarithmic vector-valued modular form and the concomitant notions of logarithmic meromorphic, holomorphic
(i.e., entire in the sense of Hecke) and cuspidal vector-valued modular forms (Subsection 2.3).

\bigskip
The Poincar\'{e} series is an indispensable device in every theory of modular (or automorphic)
forms, regardless of the level of abstraction. Naturally, then, we introduce appropriately constructed
Poincar\'{e} series to establish the existence of nontrivial logarithmic vector-valued modular forms.
Of course, in the logarithmic case we treat here the construction is of necessity more complicated, as compared with the classical  (i.e., scalar-valued) theory or the normal vector-valued case. The principal new complexity resides in the additional matrix factor $B_{\rho}$ that must be inserted in the definition
(cf. Definition 3.1) in order to achieve the desired formal transformation properties with respect to the representation $\rho$ (Subsection 3.1, following the proof of Lemma \ref{lemmacommmatrices} ).

\bigskip
It is useful to compare Definition 3.1 with the corresponding definition in the normal case
(\cite{KM2}, display (18), pp. 1352-1353).  Definition 3.1 actually defines a \emph{matrix-valued} 
Poincar\'{e} series, each column of which is a logarithmic vector-valued modular form. In fact, the same is true of our definition in the normal case, except that in the latter case we define the Poincar\'{e} series as a \emph{single column} of the matrix-valued Poincar\'{e} series. Matrix-valued modular forms 
are a very natural generalization of vector-valued modular forms. In addition to our Poincar\'{e} series, for example, the modular Wronskian \cite{M1} is the determinant of a matrix-valued modular form. The passage from vector-valued modular forms to matrix-valued modular forms is analogous to passing from a modular linear differential equation of order $p$
(loc. cit.) to an associated system of $p$ linear differential equations of order $1$.

\bigskip
Subsections 3.2 and 3.3 are devoted, respectively, to the proof of convergence of our matrix-valued Poincar\'{e} series and the determination of the general form of their logarithmic $q$-series expansions.
Our proof of convergence requires the assumption that the eigenvalues of $\rho(T)$ have absolute value 
$1$, so that the $q$-series $h_j(\tau)$ in (\ref{logform1}) again have at worst real exponents. This condition will be implicitly assumed in the remainder of the Introduction.

\bigskip
The remainder of the paper is devoted to applications. In Subsection
3.4 we give some consequences of an algebraic nature. We show (Theorem \ref{thmfreeMmod}) that if  $\rho$ has dimension $p$, the graded space
$\mathcal{H}(\rho)$ of all holomorphic vector-valued modular forms associated to $\rho$
is a \emph{free module of rank $p$} over the algebra $\mathcal{M}$ of (scalar) holomorphic modular forms on $\Gamma$. This generalizes the corresponding Theorem proved in \cite{MM} in the normal case. In fact, the proof in \cite{MM} was organized with just such a generalization in mind. The only additional input that is required is the existence of \emph{some} nonzero holomorphic vector-valued modular form associated with $\rho$, and this is an easy consequence of the existence of a nonzero meromorphic Poincar\'{e} series. A consequence of the free module Theorem is Theorem \ref{MLDEcor}, which implies
 that  if $F$ is a logarithmic vector-valued modular form $F$  then there is a \emph{canonical}
 modular linear differential equation  whose solution space is spanned by the component functions of $F$.

\bigskip
In Section 4 we derive polynomial growth estimates for the coefficients of entire
and cuspidal logarithmic vector-valued modular forms associated to $\rho$. The method here follows the approach used in \cite{KM1}, itself an extension of Hecke's venerable method for estimating the Fourier coefficients of scalar modular/automorphic forms, combined with a simple new estimate (Proposition \ref{propnlest}) that we use in Subsection 3.2 to prove convergence of our Poincar\'{e} series.

\bigskip
The occurrence of $q$-expansions of the form (\ref{logform1}) is well known in
rational and logarithmic conformal field theory. Indeed, much of the motivation for the present work originates from a need to develop a systematic theory of vector-valued modular forms wide enough in scope to cover possible applications in such field theories. By results in \cite{DLM} and \cite{M}, the eigenvalues of $\rho(T)$ for the representations that arise in rational and logarithmic conformal field theory are indeed of absolute value $1$ (in fact, they are roots of unity). Thus this assumption is natural from the perspective of conformal field theory. Our earlier results \cite{KM1} on polynomial estimates for Fourier coefficients of entire vector-valued modular forms in the normal case have found a number of applications to the theory of rational vertex operator algebras, and we expect that the extension to the logarithmic case that we prove here will be useful in the study of $C_2$-cofinite vertex operator algebras, which constitute the algebraic underpinning of logarithmic field theory.

\bigskip
Other properties of logarithmic vector-valued modular forms are also of interest, from both a foundational and applied perspective. These include
a Petersson pairing, generation of the space of cusp-forms by Poincar\'{e} series, existence of a natural boundary for the component functions, and explicit formulas (in terms of Bessel functions and Kloosterman sums) for the Fourier coefficients of Poincar\'{e} series. This program was carried through in the normal case in \cite{KM2}. We expect that the more general logarithmic case will yield a similarly rich harvest, but one must expect more complications. For example, there are logarithmic vector-valued modular forms with nonconstant component functions that may be extended to the whole of the complex plane, so that the usual natural boundary result is false \emph{per se}. Furthermore, our preliminary calculations indicate that the explicit formulas exhibit genuinely new features. We hope to return to these questions in the future.

 \section{Logarithmic vector-valued modular forms}
 \subsection{Unrestricted vector-valued modular forms}\label{UVVMF}
We start with some notation that will be used throughout.
 The  \emph{modular group} is
 \begin{eqnarray*}
\Gamma = \left\{ \left(\begin{array}{cc}a&b \\c&d\end{array}\right)\ | \ a, b, c, d \in \mathbb{Z}, \  ad-bc=1 \right\}.
\end{eqnarray*} 
It is generated by the matrices
\begin{eqnarray}\label{STgens}
S = \left(\begin{array}{cc}0 & -1 \\ 1 & 0\end{array}\right), \ T = \left(\begin{array}{cc}1&1 \\0 & 1\end{array}\right). 
\end{eqnarray}
 The complex upper half-plane is 
 \begin{eqnarray*}
\frak{H} = \{ \tau \in \mathbb{C} \ | \ \Im(\tau) > 0 \}.
\end{eqnarray*}
 There is a standard left action $\Gamma \times \frak{H} \rightarrow \frak{H}$ given by
 M\"{o}bius transformations:
 \begin{eqnarray*}
\left(  \left(\begin{array}{cc}a&b \\c&d\end{array}\right), \tau \right) \mapsto \frac{a\tau+b}{c\tau + d}.
\end{eqnarray*}
 Let $\frak{F}$ be the space of holomorphic functions in $\frak{H}$.
 There is a standard $1$-cocycle $j: \Gamma \rightarrow \frak{F}$ defined by
 \begin{eqnarray*}
j(\gamma, \tau) = j(\gamma)(\tau) = c\tau+d,  \ \ \ \gamma =  \left(\begin{array}{cc}a&b \\c&d\end{array}\right).
\end{eqnarray*}

\medskip
$\rho: \Gamma \rightarrow GL(p, \mathbb{C})$ will always denote a $p$-dimensional matrix representation of $\Gamma.$
An \emph{unrestricted vector-valued modular form of weight $k$ with respect to $\rho \  (k \in \mathbb{Z})$} is a 
holomorphic function $F: \frak{H} \rightarrow\mathbb{C}^p$ satisfying
\begin{eqnarray*}
\rho(\gamma)F(\tau) = F|_k \gamma (\tau), \ \ \gamma \in \Gamma,
\end{eqnarray*}
where the right-hand-side is the usual stroke operator
\begin{eqnarray}\label{vvdef}
F|_k \gamma (\tau) = j(\gamma, \tau)^{-k}F(\gamma \tau).
\end{eqnarray}
We could take $F(\tau)$ to be \emph{meromorphic} in $\frak{H}$, but we will not consider that more general situation here. Choosing coordinates, we can rewrite (\ref{vvdef}) in the form
\begin{eqnarray}\label{ract}
 \rho(\gamma) 
\left(\begin{array}{c}f_1(\tau) \\ \vdots \\ f_p(\tau) \end{array}\right)
= \left(\begin{array}{c}f_1|_k \gamma (\tau) \\ \vdots \\ f_p|_k (\gamma) (\tau) \end{array}\right)
\end{eqnarray}
  with each $f_j(\tau) \in \frak{F}$. We also refer to $(F, \rho)$ as an unrestricted vector-valued
  modular form.

 \subsection{Logarithmic $q$-expansions}\label{sectqexp}
In this Subsection we consider the $q$-expansions associated to unrestricted vector-valued
modular forms. 
We make use of the polynomials defined  for $k \geq 1$ by
\begin{eqnarray*}
{x \choose k} = \frac{x(x-1)  \hdots (x -k+1)}{k!},
\end{eqnarray*}
and with ${x \choose 0} = 1$ and ${x \choose k} = 0$ for $k \leq -1$.

\bigskip
We consider a finite-dimensional subspace $W \subseteq \frak{F}_k$ that is \emph{invariant} under $T$,
i.e $f(\tau +1) \in W$ whenever $f(\tau) \in W$.
We introduce the $m \times m$ matrix
\begin{eqnarray}\label{Jblock}
J_{m, \lambda} = \left(\begin{array}{cccc} \lambda &&& \\  \lambda & \ddots && \\  
&\ddots  & \ddots &  \\& & \lambda & \lambda \end{array}\right),
\end{eqnarray}
i.e. $J_{i, j}=\lambda$ for $i = j$ or $j+1$ and $J_{i, j}=0$ otherwise.  
\begin{lemma} \label{lemmaJform}There is a basis of $W$ with respect to which the 
matrix $\rho(T)$ representing
$T$ is in block diagonal form
  \begin{eqnarray}\label{Jform}
\rho(T) = \left(\begin{array}{ccc}J_{m_1, \lambda_1} &   &   \\  & \ddots &  \\ & & J_{m_t, \lambda_t}\end{array}\right).
\end{eqnarray}
\end{lemma}
\begin{pf} The existence of such a representation is basically the theory of the Jordan canonical form. The \emph{usual} Jordan canonical form is similar to the above, except that the subdiagonal of each block then consists of $1$'s rather than $\lambda$'s. The $\lambda$'s that appear in (\ref{Jform})
are the eigenvalues of
$\rho(T)$, and in particular they are nonzero on account of the invertibility of $\rho(T)$. Then it is easily checked
that (\ref{Jform}) is indeed similar to the usual Jordan canonical form, and the Lemma follows. 
$\hfill \Box$
\end{pf}

\bigskip
 We refer to (\ref{Jform})  as the \emph{modified Jordan canonical form} of $\rho(T)$,
 and  $J_{m_i, \lambda_i}$ as a \emph{modified Jordan block}.
To a certain extent at least, Lemma \ref{lemmaJform} reduces the study of the functions in $W$ to 
those associated to one of the Jordan blocks. In this case we have the following basic result.

\begin{thm}\label{thmlogqexp} Let $W \subseteq \frak{F}_k$ be a $T$-invariant 
subspace of dimension $m$.
Suppose that $W$ has an ordered basis $(g_{0}(\tau), \hdots, g_{m-1}(\tau))$
with respect to which the matrix $\rho(T)$  is a single modified Jordan block $J_{m, \lambda}$. Set $\lambda = e^{2\pi i \mu}$. Then
there are $m$ convergent $q$-expansions $h_t(\tau) = \sum_{n \in \mathbb{Z}} a_t(n)q^{n+\mu}, 0 \leq t \leq m-1,$ such that
 \begin{eqnarray}\label{polyform}
g_j(\tau) = \sum_{t=0}^j  {\tau \choose t}h_{j-t}(\tau), \ 0 \leq j \leq m-1.
\end{eqnarray}
\end{thm} 

\medskip
The case $m=1$ of the Theorem is well known. We will need it for the proof of the general case, so we state it as
\begin{lemma}\label{lemmaqexp} Let  $\lambda = e^{2 \pi i \mu}$, and suppose that $f(\tau) \in \frak{F}$ satisfies $f(\tau +1) = \lambda f(\tau)$. Then $f(\tau)$
is represented by a convergent $q$-expansion
\begin{eqnarray}\label{qexp}
f(\tau) = \sum_{n \in \mathbb{Z}} a(n)q^{n+\mu}.
\end{eqnarray}
$\hfill \Box$
\end{lemma}

Turning to the proof of the Theorem, we have
\begin{eqnarray}\label{gjrecur}
g_j(\tau +1) &=& \lambda(g_j(\tau) + g_{j-1}(\tau)), \ 0 \leq j \leq m-1,
\end{eqnarray}
where $g_{-1}(\tau) = 0$.  Set
\begin{eqnarray*}
&&h_j(\tau) =  \sum_{t=0}^j (-1)^t {\tau + t -1 \choose t}g_{j-t}(\tau), \ 0 \leq j \leq m-1.
\end{eqnarray*}
These equalities can be displayed as a system of equations. Indeed, 
\begin{eqnarray}\label{SOE}
B_m(\tau) \left(\begin{array}{c}g_{0}(\tau) \\ \vdots \\ g_{m-1}(\tau)\end{array}\right) =  \left(\begin{array}{c}h_{0}(\tau) \\ \vdots \\ h_{m-1}(\tau)\end{array}\right),
\end{eqnarray}
where $B_m(x)$ is the $m \times m$ lower triangular matrix with
\begin{eqnarray}\label{Bdef}
B_m(x)_{ij} = (-1)^{i-j}{x+i-j-1 \choose i-j} .
\end{eqnarray}
Then $B_m(x)$ is invertible and
\begin{eqnarray}\label{Bdef-1}
B_m(x)^{-1}_{ij} = {x \choose i-j}.
\end{eqnarray}

We will show that each $h_j(\tau)$ has a convergent $q$-expansion. This being the case,
(\ref{polyform}) holds and the Theorem will be proved.
 Using (\ref{gjrecur}), we have
\begin{eqnarray*}
&&h_j(\tau +1) = \lambda \sum_{t=0}^j (-1)^t 
{\tau + t  \choose t}(g_{j-t}(\tau) + g_{j-t-1}(\tau)) \\
&=&\lambda \left\{  \sum_{t=0}^{j} (-1)^t 
\left(1+ \frac{t}{\tau}   \right){\tau + t -1 \choose t}g_{j-t}(\tau)   
+ \sum_{t=0}^{j} (-1)^t 
{\tau + t  \choose t} g_{j-t-1}(\tau) \right\} \\
&=&\lambda \left\{ h_j(\tau) +  \sum_{t=0}^{j} (-1)^t 
{\tau + t -1 \choose t} \frac{t}{\tau} g_{j-t}(\tau) +  \sum_{t=0}^{j} (-1)^t 
{\tau + t  \choose t} g_{j-t-1}(\tau) \right\}. 
\end{eqnarray*}
But the sum of the second and third terms  in the braces vanishes, being equal to
\begin{eqnarray*}
&&   \sum_{t=1}^{j} (-1)^t 
{\tau + t -1 \choose t} \frac{t}{\tau} g_{j-t}(\tau) + \sum_{t=1}^{j} (-1)^{t-1} 
{\tau + t -1 \choose t -1} g_{j-t}(\tau)   \\
&=& \sum_{t=1}^{j} (-1)^{t-1} 
g_{j-t}(\tau) \left\{ {\tau +t-1 \choose t-1 }  - {\tau +t-1 \choose t }\frac{t}{\tau} \right\} = 0.
\end{eqnarray*}

\medskip
Thus we have established the identity $h_j(\tau +1) = \lambda h_j(\tau)$.
By Lemma \ref{lemmaqexp}, $h_j(\tau)$ is indeed represented by a $q$-expansion
 of the desired shape, and the proof of Theorem \ref{thmlogqexp} is complete. $\hfill \Box$
 
\bigskip
We call (\ref{polyform}) a \emph{polynomial}
$q$-expansion. The space of polynomials 
spanned by ${x \choose t}, 0 \leq t \leq m-1$ is also spanned by the powers $x^t, 0 \leq t \leq m-1$.
Since $(2\pi i \tau)^t = (\log q)^t$, it follows that in Theorem \ref{thmlogqexp} we can find a basis $\{g'_j(\tau)\}$ of $W$ such that
\begin{eqnarray}\label{logform}
g'_j(\tau) = \sum_{t=0}^j (\log q)^t h'_{j-t}(\tau)
\end{eqnarray}
with $h'_t(\tau) = \sum_{n \in \mathbb{Z}} a'_t(n)q^{n+\mu}$. 
We refer to (\ref{logform}) as a
\emph{logarithmic} $q$-expansion.

 \subsection{Logarithmic vector-valued modular forms}
We say that a function $f(\tau)$ with a $q$-expansion (\ref{qexp}) is \emph{meromorphic at infinity}
  if 
   \begin{eqnarray*}
f(\tau) = \sum_{n +\Re(\mu) \geq n_0} a(n)q^{n+\mu}.
\end{eqnarray*}
That is, the Fourier coefficients $a(n)$ \emph{vanish} for exponents $n+\mu$ whose \emph{real parts} 
are small enough. A polynomial (or logarithmic) $q$-expansion 
(\ref{polyform}) is holomorphic at infinity if each of the associated ordinary $q$-expansions
$h_{j-t}(\tau)$ is holomorphic at infinity. Similarly, $f(\tau)$ \emph{vanishes} at $\infty$ if the Fourier coefficients $a(n)$ vanish for $n+ \Re(\mu)\leq 0$;  a polynomial $q$-expansion vanishes at $\infty$ if the associated ordinary $q$-expansions vanish at $\infty$. These conditions are independent of the chosen representations.

  \medskip
  Now assume that $F(\tau) = (f_1(\tau), \hdots, f_p(\tau))^t$ is an unrestricted vector-valued modular form of weight $k$ with respect to $\rho$. It follows from (\ref{ract}) that the span $W$ of the functions $f_j(\tau)$ is a right $\Gamma$-submodule of $\frak{F}$ satisfying $f_j(\tau +1) \in W$. Choose a basis of $W$ so that $\rho(T)$ is in modified Jordan canonical form. By
  Theorem \ref{thmlogqexp} the basis of $W$ consists of functions $g_j(\tau)$ which have 
  polynomial $q$-expansions. We call $F(\tau)$, or $(F, \rho)$, a \emph{logarithmic meromorphic, holomorphic, or cuspidal vector-valued modular form}, respectively, 
  if each of the functions $g_j(\tau)$ is meromorphic, is holomorphic, or vanishes at $\infty$, respectively.
 
 \medskip
 From now on we generally drop the adjective `logarithmic' from this terminology, and say that $F(\tau)$ is \emph{semisimple} if the component functions have ordinary $q$-expansions, i.e. they are free of logarithmic terms. This holds if, and only if, $\rho(T)$ is a semisimple operator.
 
 \bigskip
 Let $\mathcal{H}(k, \rho)$ be the space of holomorphic vector-valued modular forms of weight $k$ with respect to $\rho$, with
$\mathcal{H}(\rho) = \oplus_{k \in \mathbb{Z}} \mathcal{H}(k, \rho)$
the $\mathbb{Z}$-graded space of all holomorphic vector-valued modular forms.

 \subsection{Matrix-valued modular forms}
 \emph{Matrix-valued modular forms} are a natural generalization of vector-valued modular forms. They arise naturally in several contexts, including (as we shall see) Poincar\'{e} series.
 Let $\rho: \Gamma \rightarrow GL_p(\mathbb{C})$ be a representation, and let 
$Mat_{p \times n}(\mathbb{C})$ be the space of $p \times n$ matrices. 
Let $\underline{k} = (k_1, \hdots, k_n) \in \mathbb{Z}^n.$ Consider a holomorphic map 
 $A: \frak{H} \rightarrow Mat_{p \times n}(\mathbb{C})$ satisfying
 \begin{eqnarray*}
\rho(\gamma)A(\tau) = A|_{\underline{k}} \gamma(\tau), \ \ \gamma \in \Gamma,
\end{eqnarray*}
where the right hand side is defined as
\begin{eqnarray*}
A|_{\underline{k}} \gamma(\tau) = A(\gamma \tau) J_{\underline{k}}(\gamma, \tau)^{-1}
\end{eqnarray*}
and $J$ is the matrix automorphy factor
\begin{eqnarray}\label{Jdef}
 J_{\underline{k}}(\gamma, \tau) =
 \left(\begin{array}{ccc}j(\gamma, \tau)^{k_1} & \hdots & 0 \\ \vdots & \ddots & 0 \\0 & 0 & j(\gamma, \tau)^{k_n}\end{array}\right).
\end{eqnarray}
 This defines an unrestricted matrix-valued modular form of weight $\underline{k}$ with respect to
 $\rho$.  Let $p_j: Mat_{p \times n}(\mathbb{C}) \rightarrow Mat_{p \times 1}(\mathbb{C})$ be projection onto the $j$th. column. Then $p_j \circ A$ is an unrestricted vector-valued modular form of weight $k_j$ with respect to $\rho$, and we say that $A(\tau)$ is a meromorphic, holomorphic, or cuspidal vector-valued modular form of weight $\underline{k}$ if each $p_j \circ A$ is meromorphic, holomorphic or cuspidal, respectively. Thus, a matrix-valued modular form associated to $\rho$ consists of $n$
 vector-valued modular forms of weight $k_1, \hdots, k_n$, each associated to $\rho$ with the component functions organized into the columns of a matrix.

\subsection{The nontriviality condition}\label{subsectionnontrivcond}
Let $\rho: \Gamma \rightarrow GL_p(\mathbb{C})$ be a matrix representation. Because
$S^2 = -I_2$  has order $2$, we can choose a basis of the underlying representation space such that
\begin{eqnarray*}
\rho(S^2) = \left(\begin{array}{cc}I_{p_1} & 0 \\  0& -I_{p_2}\end{array}\right).
\end{eqnarray*}
Since $S^2$ is in the \emph{center} of $\Gamma$ then $\rho(\Gamma)$ acts on the two eigenspaces of $\rho(S^2)$ and therefore the matrices
\begin{eqnarray}\label{rhodecomp}
\rho(\gamma) =  \left(\begin{array}{cc} \rho_{11}(\gamma) & 0 \\ 0 & \rho_{22}(\gamma)\end{array}\right), \ \gamma \in \Gamma,
\end{eqnarray}
are correspondingly in block diagonal form. It follows that if $(F, \rho)$ is a vector-valued modular form of weight $k$, and if we write $F(\tau) = (F_1(\tau), F_2(\tau))$ with $F_i(\tau)$ having $p_i$ components,
$i = 1, 2$, then $F_i(\tau)$ is a vector-valued modular form of weight $k$ with respect to the representation $\rho_{ii}$. More is true. The equality
$\rho(S^2)F^t(\tau) = F^t|_k S^2 (\tau)$ says that
  \begin{eqnarray*}
(F_1(\tau), -F_2(\tau)) = (-1)^k(F_1(\tau), F_2(\tau)).
\end{eqnarray*}
Assuming that $F \not= 0$, it follows that either $F_2 = 0$ and $k$ is \emph{even}, or else
$F_1 = 0$ and $k$ is \emph{odd}.
  It follows that there are natural identifications
  \begin{eqnarray}\label{Hdirsum1}
\mathcal{H}(k, \rho) &=& \left\{ \begin{array}{ll} 
                                               \mathcal{H}(k, \rho_{11}) & k \ \mbox{even},\\
                                               \mathcal{H}(k, \rho_{22}) & k \ \mbox{odd},   \notag
                                                \end{array} \right. \\
         \mathcal{H}(\rho) &=& \mathcal{H}(\rho_{11}) \oplus  \mathcal{H}(\rho_{22}).                                        
\end{eqnarray}

  \medskip
  The upshot of this discussion is that for most considerations, we may assume that
 $\rho(S^2)$ is a \emph{scalar}, i.e.
 \begin{eqnarray}\label{epsilondef}
\rho(S^2) = \epsilon I_p, \ \epsilon = \pm 1.
\end{eqnarray}
 In this case,  if $F(\tau) \in \mathcal{H}(k, \rho)$ is  nonzero then 
 \begin{eqnarray}\label{nontrivdef}
\epsilon = (-1)^k.
\end{eqnarray}
This is the \emph{nontriviality condition} in weight $k$.

\medskip
In the case of semisimple vector-valued modular forms,  it is proved in \cite{KM1} and \cite{M1} that there is an integer
$k_0$ such that $\mathcal{H}(k, \rho) = 0$ for $k<k_0$. The proof in \cite{M1} applies
to the general (logarithmic) case. Thus if $\rho$ satisfies (\ref{epsilondef}) then 
\begin{eqnarray}\label{Hdirsum2}
\mathcal{H}(\rho) = \bigoplus_{k \geq k_0} \mathcal{H}(k_0 + 2k).
\end{eqnarray}

 \section{Matrix-valued Poincar\'{e} Series}
We develop a theory of Poincar\'{e} series in order to prove \emph{existence} of nontrivial 
logarithmic vector-valued modular forms.  

\subsection{Definition and formal properties}
Fix a representation $\rho: \Gamma
\rightarrow GL(p, \mathbb{C})$. We may, and shall, assume that $\rho(T)$ is in modified Jordan canonical form 
with $t$ blocks, the $r$th. block being the $m_r \times m_r$ matrix $J_{m_r, \lambda_r}$ (\ref{Jblock}), (\ref{Jform}) and with $\lambda_r = e^{2 \pi i \mu_r}$ the associated eigenvalue of $\rho(T)$.

\medskip
We will need several more block diagonal matrices. The matrices in question will all have $t$ blocks, the $r$th block having the same size as the $r$th. block of $\rho(T)$. Set
\begin{eqnarray}\label{Brhodef}
 B_{\rho}(x) = \mbox{diag}(B_{m_1}(x), \hdots, B_{m_t}(x)),
 \end{eqnarray}
 where $B_m(x)$ is given in (\ref{Bdef}). For $(z_1, \hdots, z_t) \in \mathbb{C}^t$ let
\begin{eqnarray}\label{Lambda1def}
\Lambda_{\rho}(z_1, \hdots, z_t) = \mbox{diag}(z_1I_{m_1}, \hdots, z_tI_{m_t}).
\end{eqnarray}

\begin{defn}\label{PS}
Let $\underline{\nu} = (\nu_1, \hdots, \nu_t) \in \mathbb{Z}^t,  \underline{k} = (k_1, \hdots, k_p) \in \mathbb{Z}^p.$ The Poincar\'{e} series is defined to be
\begin{eqnarray}\label{PSdef}
 P_{\underline{k}}(\underline{\nu},  \tau) = 1/2 \sum_{M}  \rho(M)^{-1} 
\Lambda_{\rho} \left(\hdots,  e^{2 \pi i (\nu_r + \mu_r)M\tau}, \hdots \right)B_{\rho}(M\tau)^{-1}
J_{\underline{k}}(M, \tau)^{-1}, 
\end{eqnarray}
where $M$ ranges over a set of  representatives of the coset space
$\langle T \rangle \backslash \Gamma$ and $J_{\underline{k}}(M, \tau)$ is the matrix automorphy
factor (\ref{Jdef}).
\end{defn}
$B_{\rho}(\tau)^{-1}$ should be considered as an additional matrix automorphy factor.
At least formally,  $P_k(\nu, \tau)$ is a $p \times p$ \emph{matrix-valued} function.

\bigskip
We interpolate a Lemma.
\begin{lemma}\label{lemmacommmatrices} The matrices $\rho(T), \Lambda_{\rho}(z_1, \hdots, z_t)$
and $B_{\rho}(\tau) \ (\tau \in \frak{H})$ commute with each other, and satisfy
\begin{eqnarray*}
\rho(T)B_{\rho}(\tau)^{-1} = B_{\rho}(\tau +1)^{-1} \Lambda_{\rho}(\lambda_1, \hdots, \lambda_t).
\end{eqnarray*} 
\end{lemma}
\begin{pf} All of the matrices in question are block diagonal with corresponding blocks of the same size. So it suffices to show that for a given $m$ and
 $\lambda$,
 the $m \times m$ matrices $J_{m, \lambda}, zI_m$ and $B_m(\tau)$ commute and satisfy
 \begin{eqnarray}\label{morecommute}
J_{m, \lambda}B_m(\tau)^{-1} = \lambda B_m(\tau + 1)^{-1}.
\end{eqnarray}
  The $m \times m$ matrices all have the following properties: they are lower triangular and 
 the $(i, j)$-entry \emph{depends only on $i$-$j$}. It is easy to check that any two such matrices  commute. 
 
 \medskip
 As for (\ref{morecommute}), let $G(\tau)$ and $H(\tau)$ denote the column vectors of functions that occur in (\ref{SOE}), so that we can write the equation as
 \begin{eqnarray*}
B_m(\tau)G(\tau) = H(\tau).
\end{eqnarray*}
By definition of $G(\tau)$ and $H(\tau)$  (cf. Theorem \ref{thmlogqexp}) we have
\begin{eqnarray*}
J_{m, \lambda}G(\tau) &=& G(\tau +1), \\
H(\tau + 1) &=& \lambda H(\tau).
\end{eqnarray*}
Therefore, 
\begin{eqnarray*}
&&J_{m, \lambda}B_m(\tau)^{-1}H(\tau) = J_{m, \lambda}G(\tau) = G(\tau + 1) \\
&=& B_m(\tau +1)^{-1}H(\tau +1) = \lambda B_m(\tau + 1)^{-1}H(\tau).
\end{eqnarray*}
Since the components of $H(\tau)$ are linearly independent, (\ref{morecommute}) follows.
  $\hfill \Box$
\end{pf}

\bigskip
Now make the replacement $M \mapsto TM$ in a summand of (\ref{PSdef}). Using Lemma \ref{lemmacommmatrices}
we calculate that the summand maps to
\begin{eqnarray*}
&&\rho(TM)^{-1} \Lambda_{\rho} \left(\hdots,  e^{2 \pi i (\nu_r + \mu_r)TM\tau}, \hdots \right)
 B_{\rho}(TM\tau)^{-1}J_{\underline{k}}(TM, \tau)^{-1}\\
=&&\rho(M)^{-1} \rho(T)^{-1} 
\Lambda_{\rho} \left(\hdots,  e^{2 \pi i (\nu_r + \mu_r)M\tau}, \hdots \right)
\Lambda_{\rho}(\lambda_1, \hdots, \lambda_t)
B_{\rho}(M\tau +1)^{-1}J_{\underline{k}}(M, \tau)^{-1} \\
=&&\rho(M)^{-1} 
\Lambda_{\rho} \left(\hdots,  e^{2 \pi i (\nu_r + \mu_r)M\tau}, \hdots \right)
B_{\rho}(M\tau)^{-1} J_{\underline{k}}(M, \tau)^{-1}.
\end{eqnarray*}
This calculation confirms  that the sum defining $P_{\underline{k}}(\underline{\nu}, \tau)$ is \emph{independent} of the choice of coset representatives.
 We also note that
\begin{eqnarray*}
P_{\underline{k}}|_{\underline{k}} \gamma(\tau) &=& 1/2\sum_{M} \rho(M)^{-1} 
\Lambda_{\rho} \left(\hdots, e^{2 \pi i (\nu_r + \mu_r)M \gamma \tau}, \hdots \right) 
B_{\rho}(M \gamma \tau)^{-1}J_{\underline{k}}(M, \gamma \tau)^{-1} 
J_{\underline{k}}(\gamma, \tau)^{-1} \\
&=&1/2\rho(\gamma)  \sum_{M} \rho(M \gamma)^{-1} 
\Lambda_{\rho} \left(\hdots, 
 e^{2 \pi i (\nu_r + \mu_r)M \gamma \tau}, \hdots \right) B_{\rho}(M \gamma \tau)^{-1} 
 J_{\underline{k}}(M\gamma, \tau)^{-1}\\
&=& \rho(\gamma)P_{\underline{k}}(\underline{\nu}, \tau),
\end{eqnarray*}
where we used independence of coset representatives for the last equality. This confirms
that each $P_{\underline{k}}(\underline{\nu}, \tau)$ is, at least  formally, a matrix-valued modular form
 of weight $\underline{k}$ with respect to $\rho$.

\subsection{Convergence of $P_{\underline{k}}(\nu, \tau)$}\label{subsecconverge}

From now on we assume that the constants $\mu_r$ are \emph{real}, i.e. the eigenvalues $\lambda$ of 
$\rho(T)$ satisfy $|\lambda | = 1$.
With this assumption, we show in this Subsection  that the Poincar\'{e} series 
$P_{\underline{k}}(\underline{\nu}, \tau)$
is an unrestricted matrix-valued modular form for $\underline{k} \gg 0$. After the results of the previous Subsection, this amounts to the fact that $P_{\underline{k}}(\underline{\nu}, \tau)$ is holomorphic in $\frak{H}$ as long as the component weights $k_j$ of $\underline{k}$ are large enough.

\medskip
Define the vertical strip
\begin{eqnarray*}
\mathcal{S} = \{ \tau \in \frak{H} \ | \ |\Re({\tau})| \leq 1/2, \ \Im(\tau) \geq \sqrt{3}/2 \}.
\end{eqnarray*}
Notice that $\mathcal{S}$ contains the closure of the standard fundamental region 
for  $\Gamma$. We will prove
 \begin{thm}\label{thmPconv} $P_{\underline{k}}(\underline{\nu}, \tau)$ converges absolutely-uniformly 
 in $\mathcal{S}$ for $\underline{k} \gg 0$.
 \end{thm}
It is a consequence of Theorem \ref{thmPconv} and the formal transformation law for 
$P_{\underline{k}}(\underline{\nu}, \tau)$ (cf. Subsection 3.1) that $P_{\underline{k}}(\nu, \tau)$ is holomorphic throughout $\frak{H}$.

\medskip
We split off the  two terms of the Poincar\'{e} series corresponding to $\pm I$, so that
\begin{eqnarray}\label{PS1}
&& P_{\underline{k}}(\underline{\nu},  \tau) 
=  \Lambda_{\rho} \left(\hdots,  e^{2 \pi i (\nu_r + \mu_r)\tau}, \hdots \right)B_{\rho}(\tau)^{-1}+  \\
&& \ \ \ \ \ \ \ \ \ \ \ \ \ \ 1/2 \sum_{M \in \mathcal{M}^*}  \rho(M)^{-1} 
\Lambda_{\rho} \left(\hdots,  e^{2 \pi i (\nu_r + \mu_r)M\tau}, \hdots \right)B_{\rho}(M\tau)^{-1}
J_{\underline{k}}(M, \tau)^{-1}. \notag
\end{eqnarray}
Here, $\mathcal{M}^*$ is a set of representatives of the cosets 
$ \langle T \rangle \backslash \Gamma$
distinct from $\pm \langle T \rangle$, and $\pm I$ are the representative of $\pm \langle T \rangle$.
The matrices $M \in \mathcal{M}^*$ have bottom row $(c, d)$ with $c \not= 0$.
The entries of $B_{\rho}(\tau)^{-1}$ are polynomials in $\tau$  (cf. (\ref{Bdef-1}), (\ref{Brhodef})), so the first term in (\ref{PS1}) is holomorphic.

\medskip

\begin{lemma}\label{lemmaunifbd}We can choose coset representatives $M \in \mathcal{M}^*$ so that
$|M \tau|$ is \emph{uniformly bounded} in $\mathcal{S}$. That is, there is a constant $K$ such that
$|M\tau| \leq K$ for all $\tau \in \mathcal{S}$ and all $M \in \mathcal{M}^*$.
\end{lemma}
\begin{pf} Suppose that $\gamma = \left(\begin{array}{cc} a & b \\ c & d \\  \end{array}\right) \in \Gamma$ with
$c \not = 0$, and consider the $\gamma$-image $\gamma(\mathcal{S})$ of the strip. Apart from two exceptional cases (when $c= \pm 1, d=  \mp 1$), the $\gamma$-image of the circle $|\tau| = 1$  is a circle with center $b/d$ and radius at most $1$. Moreover $\gamma(\infty) = a/c$ lies on or inside the boundary of this circle. From this it is easy to see that
\begin{eqnarray*}
\gamma(\mathcal{S}) \subseteq \{ \tau \in \frak{H} \ | \ |\Re(\tau) - a/c | \leq 1 \}.
\end{eqnarray*}
Replacing $\gamma$ by $T^l \gamma$ for suitable $l$, the corresponding value of $|a/c|$ can be made less than $1$, so that
\begin{eqnarray*}
\gamma(\mathcal{S}) \subseteq \{ \tau \in \frak{H} \ | \ |\Re(\tau)| \leq 2 \}.
\end{eqnarray*}
This also holds in the exceptional cases. Therefore we may, and shall,  choose a set of coset representatives 
$\mathcal{M}^*$ so that $|\Re(M\tau)|$ is uniformly bounded
in $\mathcal{S}$ for $M \in \mathcal{M}^*$.

\medskip
On the other hand, it is easy to see that we always have $|c\tau + d|^2 \geq c^2\Im(\tau)^2$.
Since $|\Im(\tau)| \geq \sqrt{3}/2$ for $\tau \in \mathcal{S}$, it follow that
\begin{eqnarray*}
\Im(\gamma \tau) = \frac{\Im(\tau)}{|c\tau + d|^2} \leq \frac{\Im(\tau)}{c^2\Im(\tau)^2}=\frac{1}{c^2\Im(\tau)}
\leq \frac{1}{\Im(\tau)} \leq \frac{2}{\sqrt{3}},
\end{eqnarray*}
so that $|\Im(\gamma \tau)|$ is uniformly bounded in $\mathcal{S}$ for $c \not= 0$. Therefore, with our earlier choice of $\mathcal{M}^*$, it follows that
$|M\tau|$ is also uniformly bounded in $\mathcal{S}$.
This completes the proof of the Lemma. $\hfill \Box$
\end{pf}

\bigskip
Henceforth, we assume that $\mathcal{M}^*$ satisfies the conclusion of
Lemma \ref{lemmaunifbd}. 

\begin{cor}\label{corLBest}The entries of the matrices $\Lambda_{\rho} \left(\hdots, e^{2 \pi i (\nu_r + \mu_r)M \tau}, \hdots \right)$ and  $B_{\rho}(M \tau)^{-1}$ are uniformly bounded in 
$\mathcal{S}$ for $M \in \mathcal{M}^*$.
\end{cor}
\begin{pf} For $\Lambda_{\rho} \left(\hdots, e^{2 \pi i (\nu_r + \mu_r)M \tau}, \hdots \right)$ the assertion is an immediate consequence of Lemma \ref{lemmaunifbd}. As for $B_{\rho}(M \tau)^{-1}$,
we have already pointed out that it has polynomial entries, indeed the $(i, j)$-entry is 
${M\tau \choose i-j}$.
Uniform boundedness in this case is then also a consequence of Lemma \ref{lemmaunifbd}.
$\hfill \Box$
\end{pf}

\bigskip
Next we state a modification of  (\cite{E}, p. 169, display (4)) which we call  \emph{Eichler's canonical form} for elements of $\Gamma$:
\begin{lemma}\label{lemmaEprod} Each $\gamma \in \Gamma$ has a unique representation
\begin{eqnarray}\label{Eprod}
\gamma = (ST^{l_v}) \hdots (ST^{l_1})(ST^{l_0})
\end{eqnarray}
such that $(-1)^{j-1}l_j > 0$ for $1 \leq j \leq v$. $\hfill \Box$
\end{lemma}
Thus $l_1$ is positive, the $l_j$ alternate in sign for $j \geq 1$, and there is no condition on $l_0$.

\medskip
With $\gamma$ fixed for now, we  set
\begin{eqnarray}\label{Pdef}
P_0 &=& ST^{l_0},  \notag\\
P_{j+1} &=& (ST^{l_{j}})P_{j}, \ 0 \leq j \leq v-1,\notag \\
P_j &=& \left(\begin{array}{cc}a_j & b_j \\ c_j & d_j\end{array}\right), \ 0 \leq j \leq v, \\
\gamma &=&  \left(\begin{array}{cc}a & b \\ c& d\end{array}\right).
\end{eqnarray}

\begin{propn}\label{propnlest} We have
\begin{eqnarray}\label{lest}
&&|l_0l_1 \hdots l_v| \leq |d| \ \ \ \ \ \ \ \mbox{if} \ l_0 < 0; \notag \\
&&\ \  |l_1 \hdots l_v| \leq |d-c| \ \ \mbox{if} \ l_0 = 0;  \\
&&|l_0l_1 \hdots l_v| \leq |c|+ |d| \ \mbox{if} \ l_0 > 0. \notag
\end{eqnarray}
\end{propn}
\begin{pf} 

\noindent
Case A: $l_0 < 0$. We will prove by induction on $j \geq 0$ that
\begin{eqnarray}\label{indonj}
&&\ \ (i) \  |l_0l_1\hdots l_j| \leq |d_j|, \\
&&\ (ii) \   (-1)^jb_jd_j \geq 0. \notag
\end{eqnarray}
Once this is established, the case $j=v$ of (\ref{indonj})(i) proves (\ref{lest}) in Case A.
Now
\begin{eqnarray*}
P_0 = \left(\begin{array}{cc}0 & -1 \\  1 & l_0\end{array}\right), 
\end{eqnarray*}
and the case $j=0$ is clear. For the inductive step, we have
\begin{eqnarray}\label{Pj+1}
P_{j+1} = \left(\begin{array}{cc}0 & -1 \\  1 & l_{j+1} \end{array}\right)
\left(\begin{array}{cc}a_j & b_j \\ c_j& d_j\end{array}\right) 
= \left(\begin{array}{cc} -c_j & -d_j \\  a_j+l_{j+1}c_j& b_j+l_{j+1}d_j\end{array}\right).
\end{eqnarray}
Thus $(-1)^{j+1}b_{j+1}d_{j+1} = (-1)^jb_jd_j+(-1)^jl_{j+1}d_j^2 \geq 0$ where the last inequality uses induction and the inequality stated in Lemma \ref{lemmaEprod}. So (\ref{indonj})(ii) holds. 

\medskip
As for (\ref{indonj})(i), note that because $(-1)^jb_jd_j$ and $(-1)^jl_{j+1}d_j^2 $ are both nonnegative then $b_j$ and $l_{j+1}d_j$ have the \emph{same sign}. Therefore using induction again, we have $|l_0l_1 \hdots l_{j+1}|\leq |d_jl_{j+1}|\leq|b_j|+|l_{j+1}d_j|
=|b_j+l_{j+1}d_j| = |d_{j+1}|.$ This completes the proof of Case A.

\bigskip
\noindent
Case B: $l_0=0$. Notice in this case that $\gamma T^{-1} = (ST^{l_v}) \hdots (ST^{l_1})(ST^{-1})$,
which falls into Case A with $l_0=-1$. Since 
\begin{eqnarray*}
\gamma T^{-1} = \left(\begin{array}{cc}a & b \\  c & d \end{array}\right)\left(\begin{array}{cc} 1 & -1 \\  0 & 1\end{array}\right) =\left(\begin{array}{cc} a & b-a \\  c & d-c\end{array}\right)
\end{eqnarray*}
 it follows from Case A that $|l_1 \hdots l_v| \leq |d-c|$, as was to be proved.
 
 \bigskip
 \noindent
 Case C: $l_0>0.$
 We will prove by induction on $j$ that
\begin{eqnarray}\label{indonj1}
&&\ \ (i) \  |l_0l_1\hdots l_j| \leq |c_j|+|d_j|, \ j \geq 0 \\
&&\ (ii) \   (-1)^jb_jd_j, (-1)^ja_jc_j \geq 0, \ j \geq 1. \notag
\end{eqnarray}
Once again, the case $j=v$ of (\ref{indonj1})(i) proves (\ref{lest}) in Case C, and this will complete the proof of the Proposition. 
Now
\begin{eqnarray*}
P_0 = \left(\begin{array}{cc}0 & -1 \\  1 & l_0\end{array}\right), 
P_1 = \left(\begin{array}{cc} -1 & -l_0 \\  l_1& l_0l_1-1\end{array}\right).
\end{eqnarray*}
So when $j=0$,  (\ref{indonj1})(i) is clearly true, and because  $l_0, l_1>0$ we also have
\begin{eqnarray*}
-a_1c_1= l_1> 0, \ -b_1d_1= l_0(l_0l_1-1) \geq 0.
\end{eqnarray*}
So (\ref{indonj1})(ii) holds for  $j=1$. As for the inductive step, $P_{j+1}$ is as in
(\ref{Pj+1}), and the proof that $(-1)^jb_jd_j\geq 0$ is the same as in case A. Similarly
$(-1)^{j+1}a_{j+1}c_{j+1} = (-1)^jc_ja_j+(-1)^{j}l_{j+1}c_j^2 \geq 0$ is the sum of two nonnegative terms and hence is itself nonnegative, so (\ref{indonj1})(ii) holds. Finally, by an argument similar to that used in Case A, we have
$|l_0 \hdots l_{j+1}| \leq |c_j+d_j||l_{j+1}| < |l_{j+1}c_j|+|l_{j+1}d_j|+|a_j|+|b_j| =|a_j+ l_{j+1}c_j|
+|b_j+ l_{j+1}d_j|= |c_{j+1}|+|d_{j+1}|.$
The Proposition is proved. $\hfill \Box$
 \end{pf}

\bigskip 
The \emph{Eichler length}
of $\gamma$ is given by
\begin{eqnarray}\label{Ldef}
L(\gamma) = \left \{  \begin{array}{ll}
2v+2, &  \ l_0 \not= 0 \\
2v+1, &  \ l_0 = 0    \end{array}  \right.
\end{eqnarray}
By Lam\'{e}'s Theorem we have the estimate
\begin{eqnarray}\label{Lameest}
L(\gamma) \leq K(\log |c| +1)
\end{eqnarray}
with a positive constant $K$ independent of $\gamma$.

\bigskip
The \emph{norm} $||\rho(\gamma)||$, defined to be $\max_{i, j} |\rho(\gamma)_{i j}|$,  satisfies
\begin{eqnarray}\label{norm1est}
||\rho(\gamma)|| \leq ||\rho(S)||^{v+1} \prod_{j=0}^v ||\rho(T^{l_j})||. 
\end{eqnarray}

\begin{lemma}\label{lemmaTlest} Let $s$ be the maximum of the sizes $m_j$ of the Jordan blocks
$J_{m_j, \lambda_j}$ of $\rho(T)$ (\ref{Jblock}), (\ref{Jform}). There is a constant $C_s$ depending only on $s$ such that for $l\not= 0$,
\begin{eqnarray}\label{Tlest}
||\rho(T^l)|| \leq C_s|l|^{s-1}.
\end{eqnarray}
\end{lemma}
\begin{pf} We have 
\begin{eqnarray*}
J_{m, \lambda}^l = \lambda^lJ_{m, 1}^l = \lambda^l(I_m + N)^l = \lambda^l \sum_{i\geq 0} {l \choose i}N^i
\end{eqnarray*}
where $N$ is the nilpotent $m \times m$ matrix with each $(i, i-1)$-entry equal to $1 \ (i \geq 2),$
and all other entries zero. Note that $N^m = 0$ and the entries of $N^i$ for $1 \leq i<m$ are 
$1$ on the $i$th. subdiagonal and zero elsewhere. Bearing in mind that $|\lambda| = 1$, it follows that
$||J_{m, \lambda}^l||$ is majorized by the maximum of the binomial coefficients
${l \choose i}$ over the range $0 \leq i \leq m-1$. Since ${l \choose i}$ is a polynomial in $l$ of degree
$i$ then we certainly have $||J_{m, \lambda}^l|| \leq C_m|l|^{m-1}$ for a universal constant $C_m$, and since this applies to each Jordan block of $\rho(T^l)$ then the Lemma follows immediately.
$\hfill \Box$
\end{pf}

\begin{cor}\label{lemmapolyest} There are universal constants $K_3, K_4$ such that
\begin{eqnarray}\label{polyest}
||\rho(\gamma)||  \leq K_3(c^2+d^2)^{K_4}.
\end{eqnarray}
Moreover the \emph{same} estimate holds for $||\rho(\gamma^{-1})||$.
\end{cor}
\begin{pf}
From Lemma \ref{lemmaTlest} and (\ref{norm1est}) we obtain
\begin{eqnarray*}
||\rho(\gamma)|| \leq \left \{ \begin{array}{ll} 
  K_1^v\prod_{j=0}^v |l_j|^{s-1}, & l_0 \not= 0, \\
  K_1^v\prod_{j=1}^v |l_j|^{s-1}, & l_0=0,
                                    \end{array}  \right.
\end{eqnarray*}
for a constant $K_1$ depending only on $\rho$. Now use (\ref{Ldef}), (\ref{Lameest}) 
and Proposition \ref{propnlest} to see that
\begin{eqnarray*}
||\rho(\gamma)|| \leq e^{(\log K_1)K_2 \log (|c|+1)}(|c|+|d|) \leq K_3(c^2+d^2)^{K_4}. 
\end{eqnarray*}

\medskip
Concerning the second assertion of the Lemma, since
$\gamma^{-1} = (T^{-l_0}S)(T^{-l_1}) \hdots (T^{-l_v}S)$
then
\begin{eqnarray*}
||\rho(\gamma^{-1})|| \leq  ||\rho(S)||^{v+1}\prod_{j=0}^v ||\rho(T)^{-l_j}||,
\end{eqnarray*}
and (\ref{Tlest}) then holds by Lemma \ref{lemmaTlest}. The rest of the proof is the same as the previous case, so that we indeed obtain the estimate (\ref{polyest}) for $\gamma^{-1}$ as well as $\gamma$.
$\hfill \Box$
\end{pf}

\bigskip
We can now prove Theorem \ref{thmPconv}. Let $P^*_{\underline{k}}(\underline{\nu}, \tau)$ denote the infinite
sum 
in (\ref{PS1}). We have
by Corollary \ref{corLBest} and Corollary \ref{lemmapolyest}   that
\begin{eqnarray*}
&&||P^*_{\underline{k}}(\underline{\nu},  \tau)|| \leq \sum_{M \in \mathcal{M}^*}  ||\rho(M)^{-1} ||
||\Lambda_{\rho} \left(\hdots,  e^{2 \pi i (\nu_r + \mu_r)M\tau}, \hdots \right)|| ||B_{\rho}(M\tau)^{-1}||
||J_{\underline{k}}(M, \tau)^{-1}|| \\
&& \hspace{2cm} \leq K_5 \sum_{(c, d)=1}(c^2+d^2)^{K_4} ||J_{\underline{k}}(M, \tau)^{-1}||,
\end{eqnarray*}
with constants $K_4, K_5$ that depend only on $\rho$. We also know (\cite{KM1}, display (13))
that 
\begin{eqnarray}\label{KM13}
c^2+d^2 \leq K_6|c\tau + d|^2
\end{eqnarray}
for a universal constant $K_6$.
Because of the nature of
the matrix automorphy factor $J_{\underline{k}}$ (\ref{Jdef}), it follows from the previous two
displays that if
the \emph{minimum} of the weights $k_i$ in $\underline{k} = (k_1, \hdots, k_p)$
is \emph{large enough}, then
\begin{eqnarray*}
&&||P^*_{\underline{k}}(\underline{\nu},  \tau)|| 
 \leq K_7 \sum_{(c, d)=1} (c\tau + d)^{-k},
\end{eqnarray*}
with $k > 2$. It is well-known that this series converges absolutely-unformly in $\mathcal{S}$,
so the same is true for $P^*_{\underline{k}}(\underline{\nu}, \tau)$.
This completes the proof of Theorem \ref{thmPconv}. $\hfill \Box$

\subsection{$q$-expansions of the component functions}
We now assume (cf. the discussion in Subsection \ref{subsectionnontrivcond}) that
(\ref{epsilondef}) holds. 
Consider the substitution $M \mapsto -M$ in the expression for $P_{\underline{k}}(\underline{\nu}, \tau)$.
Because the sum is independent of the order of the terms,   (\ref{epsilondef}) 
implies that  
\begin{eqnarray*}
 P_{\underline{k}}(\underline{\nu}, \tau)=P_{\underline{k}}( \underline{\nu}, \tau) \Lambda_{\rho}(\epsilon(-1)^{k_1}, \hdots, \epsilon(-1)^{k_t}).
\end{eqnarray*}
If the nontriviality condition (\ref{nontrivdef}) holds in weight $k_j$ then $\epsilon(-1)^{k_j} = 1$ and the
$j$th. column of $P_{\underline{k}}(\underline{\nu}, \tau)$ is unchanged. If the nontriviality condition does
\emph{not} hold then the $j$th. column is zero and as such it too is unchanged. 
We conclude that 
\begin{eqnarray}\label{newPS}
P_{\underline{k}}(\nu, \tau) =  \sum_{M}  \rho(M)^{-1} 
\Lambda_{\rho} \left(\hdots,  e^{2 \pi i (\nu_r + \mu_r)M\tau}, \hdots \right)B_{\rho}(M\tau)^{-1}
J_{\underline{k}}(M, \tau)^{-1},
\end{eqnarray}
where the matrices $M$ now range over an arbitrary set of coset representatives of
$\pm \langle T \rangle \backslash \Gamma$.

\medskip
We will show that $P_{\underline{k}}(\underline{\nu},  \tau)$ is a meromorphic
vector-valued modular form for 
$\underline{k} \gg 0$. We have already proved that it is an unrestricted vector-valued modular form, so that the component functions  that occur in the matrix representation
\begin{eqnarray*}
P_{\underline{k}}(\underline{\nu},  \tau) = ( P_{mn}(\tau))
\end{eqnarray*}
have polynomial $q$-expansions (\ref{polyform}) by Theorem \ref{thmlogqexp}. It remains to show that these $q$-expansions are meromorphic at infinity if the weights are large enough.
Note that because of our assumption that the constants $\mu_r$ are real, the $q$-expansions in question have only real powers of $q$. 

\medskip
To describe $P_{mn}(\tau)$, let us assume that the nontriviality condition holds in weight $k_n$.
Let  $r$ be such that the $m$th. row of $P(\tau)$ falls into the $r$th. Jordan block, and let $M_r = m_1 + \hdots + m_{r}$. Then we have
\begin{eqnarray*}
M_{r-1} < m \leq M_r.
\end{eqnarray*}
Now take $I$ as the coset representative of $\pm\langle T \rangle \backslash \Gamma$ and set 
$\mathcal{M} = \mathcal{M}^* \cup \{ I \}$. From (\ref{newPS}) we have
\begin{eqnarray*}
&&\hspace{6cm}P(\tau)_{mn} \\
&=&  \sum_{l= 1}^p \sum_{M \in \mathcal{M}} \Lambda_{\rho}( \hdots, e^{2 \pi i (\nu_s + \mu_s)M\tau}, \hdots)_{mm}
 \rho(M^{-1})_{ml} B_{\rho}(M\tau)_{ln}^{-1} j(M, \tau)^{-k_n} \\
 &=&e^{2\pi i(\nu_r + \mu_r)\tau} B_{\rho}(\tau)_{mn}^{-1}+  \sum_{l= 1}^p \left\{ \sum_{M \in \mathcal{M}^*} e^{2 \pi i (\nu_r + \mu_r)M\tau}
 \rho(M^{-1})_{ml} B_{\rho}(M\tau)^{-1}_{ln} j(M, \tau)^{-k_n} \right\}.
\end{eqnarray*}

Because of absolute-uniform convergence in the strip $\mathcal{S}$,
$\lim_{\tau \rightarrow i \infty}$ may be taken inside the summations. By Lemma \ref{lemmaunifbd}
and Corollary \ref{corLBest} we find that
\begin{eqnarray*}
\lim_{\tau \rightarrow i \infty } \left\{ P(\tau)_{mn}-e^{2\pi i(\nu_r+\mu_r)\tau}B_{\rho}(\tau)_{mn}^{-1}
\right\} = 0
\end{eqnarray*}
($k_n>0$). It follows that the polynomial $q$-expansion of $P(\tau)_{mn}-e^{2\pi i(\nu_r+\mu_r)}B_{\rho}(\tau)^{-1}$ 
can have only \emph{positive} powers of $q$, so that 
\begin{eqnarray}\label{Pmnlogqexp}
P(\tau)_{mn} &=& d_{mn} {\tau \choose m-n}q^{\nu_r + \mu_r} +
\sum_{u=0}^{m-M_{r-1}-1}{\tau \choose u} \sum_{l+\mu_r > 0}  \hat{a}_{unr}(l)q^{l+\mu_r},   \\
d_{mn} &=& \left \{ \begin{array}{ll}
                                   1, & \ M_{r-1}<n \leq m, \\
                                   0, & \ \mbox{otherwise}.
                                    \end{array}
                                     \right. \notag
\end{eqnarray}

\medskip
Notice that the diagonal terms have polynomial $q$-expansions
\begin{eqnarray*}
P(\tau)_{mm} = q^{\nu_r + \mu_r} + \mbox{regular terms}.
\end{eqnarray*}
In particular,  if $\nu_r+\mu_r<0$ then there is a pole at $i\infty$ and if
$\nu_r + \mu_r = 0$ the constant term is $1$. So in both cases $P_{mm}(\tau)$ is \emph{nonzero}.
 We have  established the following.
 \begin{thm}\label{thmPSexist} Suppose that $\rho$ satisfies $\rho(S^2) = \epsilon I_p$. Then 
 $P_{\underline{k}}(\underline{\nu}, \tau)$
 is a meromorphic matrix-valued modular form of weight $\underline{k}$ for all
 $\underline{k} \gg 0$.  If the nontriviality condition holds in all weights $k_n$ then one of the following holds: \\
 (a) \  $\nu_r + \mu_r > 0$ for all $r$ and $P_{\underline{k}}(\underline{\nu}, \tau)$ is a cuspidal matrix-valued modular form, possibly zero. \\
 (b) \ $\nu_r + \mu_r \geq 0$ for all $r, \nu_r + \mu_r = 0$ for \emph{some} $r$, and
 $P_{\underline{k}}(\underline{\nu}, \tau)$ is a \emph{nonzero}, holomorphic matrix-valued modular form
 of weight $\underline{k}$. \\
 (c)  \ $\nu_r + \mu_r < 0$ for \emph{some} $r$ and
 $P_{\underline{k}}(\underline{\nu}, \tau)$ is a \emph{nonzero} meromorphic matrix-valued modular form of weight $\underline{k}$. 

\noindent
 If the nontriviality condition is \emph{not} satisfied in weight $k_n$, then
 the $n$th. column of $P_{\underline{k}}(\underline{\nu}, \tau)$ vanishes identically. $\hfill \Box$
  \end{thm}

\subsection{Further consequences}

We record a consequence of the nature of the $q$-expansions (\ref{Pmnlogqexp}).
\begin{thm}\label{thmessentialexist} Suppose that $\rho(S^2) = \epsilon I_p$. For large enough weight $k$ there is $F(\tau) \in \mathcal{H}(k, \rho)$ such that the component functions of $F(\tau)$ are \emph{linearly independent}.
\end{thm}
\begin{pf} Let $\underline{k} = (k, \hdots, k)$ have constant weight $k$, and choose $k$ large enough so that $P(\tau) = P_{\underline{k}}(\underline{\nu}, \tau)$ is holomorphic throughout $\frak{H}$. This holds for any choice of $\underline{\nu}$. We may, and shall, also assume that the nontriviality condition in weight $k$ is satisfied. 

\medskip
Now choose $\underline{\nu}$ so that the exponents $\nu_r + \mu_r$ are negative and pairwise distinct for each $r$ in the range $1\leq r \leq t$. Consider the resulting $t$ vector-valued modular forms
$p_{M_r} \circ P(\tau) =P(\tau)_{M_r}$ where $p_{M_r}$ is projection onto the $M_r$th column (cf. Subsection 2.4). By (\ref{Pmnlogqexp}) we see that the component functions of
$P(\tau)_{M_r}$ are holomorphic outside of the $r$th. block, and in the $r$th. block they have $q$-expansions $q^{\nu_r + \mu_r}{\tau \choose u}+ \hdots, \ 0 \leq u \leq m_r-1$. Clearly, then, these functions are linearly independent. 

\medskip
Consider the vector-valued modular form
\begin{eqnarray*}
F(\tau) = \Delta(\tau)^v \sum_{r=1}^t P(\tau)_{M_r},
\end{eqnarray*}
with $v$ an integer satisfying  $v +\nu_r + \mu_r \geq 0, \ 1 \leq r \leq t$. Since the $\nu_r + \mu_r$ are pairwise distinct, it follows from the discussion of the preceding paragraph that  $\sum_r P(\tau)_{M_r}$ has linearly independent component functions. The choice of $v$ ensures that $F(\tau)$
is holomorphic at $i \infty$ and it also has linearly independent component functions. Since $F(\tau)$ is a vector-valued modular form associated with the same representation $\rho$, the Theorem follows. 
$\hfill \Box$
\end{pf}

\bigskip
As discussed in Subsection \ref{subsectionnontrivcond}, $\rho$ is equivalent to the direct sum $\rho_{1} \oplus \rho_{-1}$ of a pair of representations
$\rho_{\epsilon}$ of $\Gamma$ with the property that $\rho_{\epsilon}(S^2) = \epsilon I, \ \epsilon = \pm 1$.
From (\ref{Hdirsum1})-(\ref{Hdirsum2}) it follows that there is a natural identification
\begin{eqnarray}\label{natid}
\mathcal{H}(\rho) = \mathcal{H}(\rho_1) \oplus \mathcal{H}(\rho_{-1}),
\end{eqnarray}
with
\begin{eqnarray*}
\mathcal{H}(\rho_1) &=& \bigoplus_{k \ even} \mathcal{H}(k, \rho_1),\\
\mathcal{H}(\rho_{-1}) &=& \bigoplus_{k \ odd} \mathcal{H}(k, \rho_{-1}).
\end{eqnarray*}
In other words, $\mathcal{H}(\rho_1)$ and $\mathcal{H}(\rho_{-1})$ are the even and odd parts respectively of
$\mathcal{H}(\rho)$.

\begin{cor}\label{corvvexist} For any representation $\rho: \Gamma \rightarrow GL_p(\mathbb{C})$,
there is a \emph{nonzero} holomorphic vector-valued modular form 
$F(\tau) \in \mathcal{H}(k, \rho)$ for \emph{large enough} weight $k$.
\end{cor}
\begin{pf} If $\rho(S^2) = \pm I_p$ then the Corollary follows immediately from  
Theorem \ref{thmessentialexist}.
The general result is then a consequence of the preceding comments. $\hfill \Box$.
\end{pf}

\bigskip
Let $\mathcal{M} = \oplus_{k \geq 0}\mathcal{M}_k =  \mathbb{C}[Q, R]$ be the weighted
polynomial algebra of holomorphic modular forms of level $1$ on $\Gamma$, where $Q=E_4(\tau), R = E_6(\tau)$. As in \cite{M1}, 
\begin{eqnarray*}
\mathcal{R} = \mathcal{M}[d]
\end{eqnarray*}
is the ring of differential operators obtained by adjoining to $\mathcal{M}$ an element $d$ satisfying
\begin{eqnarray*}
df - fd = D(f), \  f \in \mathcal{M},
\end{eqnarray*}
where $D$ is the modular derivative
\begin{eqnarray}\label{Dact}
Df = D_kf = (\theta + kP)f \ \ \  (f \in \mathcal{M}_k).
\end{eqnarray}
Here, $\theta = qd/dq$ and $P = -1/12 + 2\sum_{n \geq 1}\sigma_1(n)q^n$ is the weight $2$
quasimodular Eisenstein series, normalized as indicated.

\medskip
$\mathcal{R}$ is a $2\mathbb{Z}$-graded algebra ($d$ has degree $2$), and $\mathcal{H}(\rho)$
is a $\mathbb{Z}$-graded $\mathcal{R}$-module in which $f \in \mathcal{M}$ acts as a multiplication operator and $d$ acts on $F \in \mathcal{H}(\rho)$ by its action on components of $F$ given by (\ref{Dact}). In particular, it follows that $\mathcal{R}$ operates on the even and odd parts of
$\mathcal{H}(\rho)$, so that the identification (\ref{natid}) is one of $\mathcal{R}$-modules.

\begin{thm}\label{thmfreeMmod}($\mbox{Free module Theorem}$) $\mathcal{H}(\rho)$ is a \emph{free} $\mathcal{M}$-module of \emph{rank $p$}.
\end{thm}

This means that there are $p$ weights $k_1, \hdots, k_p$ and $p$ vector-valued modular forms
$F_j(\tau) \in \mathcal{H}(k_j, \rho), \ 1 \leq j \leq p$, such that every $F(\tau) \in \mathcal{H}(k, \rho)$
has a \emph{unique} expression in the form
\begin{eqnarray*}
F(\tau) = \sum_{j=1}^p f_j(\tau) F_j(\tau), \ f_j(\tau) \in \mathcal{M}_{k-k_j}.
\end{eqnarray*}
It is an immediate consequence of this result that the Hilbert-Poincar\'{e} series for $\mathcal{H}(\rho)$
is a  rational function:
\begin{eqnarray*}
\sum_{k \geq k_0} \dim \mathcal{H}(k, \rho)t^k = \frac{\sum_{j=1}^p t^{k_j}}{(1-t^4)(1-t^6)}.
\end{eqnarray*}

\bigskip
With Corollary \ref{corvvexist} available, the remaining details of the proof of Theorem \ref{thmfreeMmod} are essentially identical to that of the semisimple case given in \cite{MM} and involve mainly arguments from commutative algebra.  In the few places where the nature of 
 the component functions of  vector-valued modular forms is relevant, the argument is the same whether the $q$-expansions are ordinary or logarithmic. We forgo further discussion.

\bigskip
We give an application of the free module Theorem. Let $F(\tau) \in \mathcal{H}(k, \rho)$. If the
elements $F, DF, \hdots, D^pF$ are linearly independent over $\mathcal{M}$ then they span a free
$\mathcal{M}$-submodule of $\mathcal{H}(\rho)$ of rank $p+1$. Since $\mathcal{H}(\rho)$ has rank $p$, this is not possible. Therefore, $F(\tau)$ satisfies an equality of the form
\begin{eqnarray}\label{FMLDE}
\left(g_0(\tau)D^p_k + g_1(\tau)D_k^{p-1} + \hdots + g_p(\tau)\right)F = 0.
\end{eqnarray}
where $g_0 \in \mathcal{M}_l$ for some weight $l$ and $g_j(\tau) \in \mathcal{M}_{l+2j}$. We may think of (\ref{FMLDE}) as a modular linear differential equation (MLDE) \cite{M1} of order at most $p$, in which case the component functions of $F(\tau)$ are solutions. Now suppose that
the component functions  are \emph{linearly independent}.  Since they are solutions of any MLDE satisfied by $F(\tau)$, the solution space must have dimension at least $p$,
and therefore the order of the MLDE must itself be at least $p$. We have therefore shown that if $F(\tau)
\in \mathcal{H}(k, \rho)$ has linearly independent component functions, it satisfies an MLDE of order 
$p$ and none of order less than $p$.

\bigskip
Continuing with the assumption that the component functions of $F(\tau)$ are linearly independent,
let $I \subseteq \mathcal{M}$ be the set of \emph{all} leading coefficients $g_0(\tau)$ that occur
in order $p$ MLDE's (\ref{FMLDE}) satisfied by $F$. Taking account of the trivial case when all coefficients
$g_j(\tau)$ vanish, we see easily that $I$ is a graded ideal. Moreover, our previous comments show that $I \not= 0$. We will show that $I$ contains  a \emph{unique} nonzero modular form $g(\tau)$ of least weight, normalized so that the leading coefficient of its $q$-expansion is $1$, and  that $I = g(\tau)\mathcal{M}$.

\bigskip
For nonzero $h_0(\tau) \in I$ of weight $k$, we let
\begin{eqnarray*}
L_h = h_0(\tau)D^p + h_1(\tau)D^{p-1} + \hdots + h_p(\tau)
\end{eqnarray*}
be the unique order $p$ differential operator in $\mathcal{R}$ with leading coefficient $h_0(\tau)$ and satisfying $L_hF = 0$. Let $g_0(\tau)$ be any nonzero element in $I$ of least weight, say $m$. Then we have
$L_gF = L_hF = 0$, and therefore also
\begin{eqnarray*}
(g_0L_h - h_0L_g)F=0.
\end{eqnarray*}
The differential operator in the last display has order at most $p-1$, and therefore (by our earlier remarks) must vanish identically. It follows that for all indices $j$ we have
\begin{eqnarray}\label{ghcompare}
g_0h_j  = h_0g_j.
\end{eqnarray}

\bigskip
Suppose that the order of vanishing of $g_0(\tau)$ at $\infty$ is \emph{greater} than that of $h_0(\tau)$.
By (\ref{ghcompare}) it follows that all $g_j(\tau)$ vanish to order at least $1$ at $\infty$, i.e. each
$g_j(\tau)$ is divisible by the discriminant $\Delta(\tau)$ in $\mathcal{M}$. But then 
$ L_gF = \Delta(\tau)L_{g'}F = 0$, whence $L_{g'}F=0$ for some nonzero $g'(\tau) \in \mathcal{M}_{m-12}$. Then $g'(\tau) \in I$, and this contradicts the minimality of the weight of $g$. Thus we have shown that the order of vanishing of $g_0(\tau)$ at $\infty$ is minimal among nonzero elements in $I$, and that this assertion holds for \emph{any} nonzero element of least weight in $I$.

\bigskip
If there are two linearly independent elements $a(\tau), b(\tau)$, say, of least weight in $I$ then some linearly combination of them
 vanishes at $\infty$ to an order that \emph{exceeds} that of at least one of $a(\tau)$ and $b(\tau)$. By the last paragraph this cannot occur, and we conclude that, up to scalars, $g_0(\tau)$ is the unique nonzero element in $I$ of least weight. 

\bigskip
We use similar arguments to show that $g_0(\tau)$ generates $I$. If not, choose an element
$h_0 \in I$ of least weight $n$, say, subject to $h_0(\tau) \notin g_0(\tau)\mathcal{M}$. If $h_0(\tau)$ has greater order of vanishing at $\infty$ than $g_0(\tau)$, (\ref{ghcompare}) and a previous argument shows that
every $h_j(\tau)$ is divisible by $\Delta(\tau)$. Then as before, $h_0(\tau) = \Delta(\tau)h'_0(\tau)$ with
$h'(\tau) \in I$. By minimality of the weight of $h_0(\tau)$ we get $h'(\tau) \in g_0(\tau)\mathcal{M}$, and therefore also $h_0(\tau) \in g_0(\tau)\mathcal{M}$, contradiction. Therefore, every element of 
 weight $n$ in $I \setminus g_0(\tau)\mathcal{M}$ has the same order of vanishing at $\infty$ as $g_0(\tau)$. This again implies the unicity of $h_0(\tau)$ up to scalars. 

\bigskip
If $n - m \geq 4$  then $h_0(\tau) + E_{n-m}(\tau)g_0(\tau)$ has weight $n$ and lies in
$I \setminus g_0(\tau)\mathcal{M}$. (Here, $E_k(\tau)$ is the usual weight $k$ Eisenstein series.)
Thus $h_0(\tau)$ is a scalar multiple of $h_0(\tau) + E_{n-m}(\tau)g_0(\tau)$ and therefore lies in
$g_0(\tau)\mathcal{M}$, contradiction.
Therefore,  $n-m = 2$. 
In this case we consider $h'(\tau) = E_4(\tau)h_0(\tau) - \beta E_6(\tau)g_0(\tau)$ and
$h''(\tau) = E_6(\tau)h_0(\tau) - \gamma E_4^2(\tau)g_0(\tau)$ for scalars $\beta, \gamma$ chosen in each case so that the order of vanishing at $\infty$ is greater than that of $g_0(\tau)$. A previous argument shows that $L_{h'}F = \Delta L_{h'_1}F=0$ for some $h_1'(\tau)$ of weight $n+4-12=m-6$.
Since $h_1'(\tau) \in I$ has weight less than $m$ then $h_1'(\tau) = 0$, so that
 $E_4(\tau)h_0(\tau) = \beta E_6(\tau)g_0(\tau)$. The same reasoning applied to $h''(\tau)$ also shows that $E_6(\tau)h_0(\tau) = \gamma E_4^2(\tau)g_0(\tau)$. From these equalities we deduce that 
$g_0(\tau)(\gamma E_4^3(\tau) - \beta E_6^2(\tau)) = 0$. This can only happen if $\beta = \gamma = 0$, whence $E_4(\tau)h_0(\tau)=0$. This is impossible since $h_0(\tau)$ is nonzero, and we have contradicted the assumed existence of $h_0(\tau)$. To summarize, we have established

\begin{thm}\label{MLDEcor} Suppose that $F(\tau) \in \mathcal{H}(k, \rho)$ has linearly independent component functions. Then the component functions are a basis of the solution space of a modular linear differential equation
\begin{eqnarray}\label{MLDE}
\left(g_0(\tau)D^p_k + g_1(\tau)D_k^{p-1} + \hdots + g_p(\tau)\right)f = 0
\end{eqnarray}
where $g_j(\tau) \in \mathcal{M}_{l+2j }, \ 0 \leq j \leq p,$ for some $l \geq 0$. The set of leading coefficients $g_0(\tau)$ that can occur in (\ref{MLDE}) is a (nonzero) principal graded ideal 
$I \subseteq \mathcal{M}$ generated by the \emph{unique} normalized modular form $g(\tau)$
of least weight in $I$. $\hfill \Box$
\end{thm}

\bigskip
If the condition that the component functions of $F(\tau)$ are linearly independent is 
\emph{not} met, one can replace $\rho$ by the representation $\rho'$ of $\Gamma$ furnished by the span of the component functions. Then the Theorem applies to $\rho'$. In this way, we see that to any logarithmic vector-valued modular form we can associate an MLDE in a canonical way: it is the MLDE of least order and with normalized leading coefficient of least weight whose solution space is spanned by the component functions of $F(\tau)$.

\bigskip
We can alternatively couch these results in terms of annihilators in the ring of differential operators
$\mathcal{R}$. For example, we have
\begin{cor} Let $F \in \mathcal{H}(k, \rho)$. Then the annihilator Ann$_{\mathcal{R}}(F)$ is a
\emph{cyclic} $\mathcal{R}$-module. $\hfill \Box$
\end{cor}

\section{Polynomial estimates of Fourier coefficients}
 Let $F(\tau) \in \mathcal{H}(k, \rho)$ be a logarithmic, holomorphic vector-valued modular form of weight $k$. We are going to show that the Fourier coefficients of $F(\tau)$
 satisfy a polynomial growth condition
for $n \rightarrow \infty$. Let $F(\tau) = (f_1(\tau), \hdots, f_p(\tau))^t$. We know by Theorem \ref{thmlogqexp}
that  there are $m_j$ $q$-expansions
$h_l(\tau) = \sum_{n+\mu_j \geq 0}a_{jl}(n)q^{n+\mu_j}, \ 0 \leq l \leq m_j-1$
such that the components of $F(\tau)$ corresponding to the $j$th Jordan block
are $(f_{m_j-1}(\tau), \hdots, f_0(\tau))^t$ with 
\begin{eqnarray*}
f_l(\tau) = \sum_{u=0}^l  {\tau \choose u} h_{l-u}(\tau), \ 0 \leq l \leq m_j-1.
\end{eqnarray*}
Here we have relabelled the components in the $j$th block for notational convenience.

\medskip
The proof is similar to the case  treated in \cite{KM1}, but with an additional complication due to the fact that we are  dealing  with polynomial $q$-expansions rather than ordinary $q$-expansions. To deal with this we make use of the
estimates  that we have obtained in Subsection \ref{subsecconverge}.
We continue to assume that the eigenvalues of $\rho(T)$ are of absolute value $1$. We will sometimes  drop the subscript $j$ from the notation when it is convenient.

\medskip
We write $\tau = x+iy$ for $\tau \in \frak{H}$ and let $\frak{R}$ be the usual fundamental region for
$\Gamma$. Write $z = u + iv$ for $z \in \overline{\frak{R}}$.
Choose a real number $\sigma > 0$ to be fixed later, and set
\begin{eqnarray*}
g_l(\tau) = y^{\sigma}|f_l(\tau)|.
\end{eqnarray*}
Because $F(\tau)$ is holomorphic, $a_l(n)=0$ unless $n+\mu \geq 0$. It follows that there is a constant $K_1$ such that
\begin{eqnarray}\label{gbound}
g_l(z) \leq K_1v^{\delta(\sigma+1)}, \ 1 \leq l \leq p, \ z \in \overline{\frak{R}},
\end{eqnarray}
where $\delta = 0$ if $F(\tau)$ is a \emph{cusp-form}, and is $1$ otherwise.

\medskip
Choose $\gamma = \left(\begin{array}{cc}a&b \\c&d\end{array}\right) \in \Gamma$, set $\tau = \gamma z$, and write $\gamma$  in Eichler canonical form (\ref{Eprod}).
We now argue just as in \cite{KM1} pp. 121-122. Thus
\begin{eqnarray*}
g_l(\tau) &=& g_l(\gamma z) =(v |cz + d|^{-2})^{\sigma}|f_l(\gamma z)| \\
&=& v^{\sigma} |cz+d|^{k-2\sigma}|f_l|_k \gamma (z)| \\
&=& v^{\sigma} |cz+d|^{k-2\sigma}|(\rho(\gamma) f (z))_l| \\
&=&  v^{\sigma} |cz+d|^{k-2\sigma}|\sum_{m=1}^p \rho(\gamma)_{lm}f^m(z)|.
\end{eqnarray*}
 Using (\ref{gbound}), Lemma \ref{lemmapolyest}, 
and (\ref{KM13}), we obtain
\begin{eqnarray*}
&&g_l(\tau) \leq K_1v^{\delta(\sigma + 1)} |cz+d|^{k-2\sigma}\sum_{m=1}^p |\rho(\gamma)_{lm}| \\
&&\ \ \ \  \ \ \ \leq K_2v^{\delta(\sigma + 1)} |cz+d|^{k-2\sigma} |c^2+d^2|^{K_4} \\
&&\ \ \ \  \ \ \ \leq K_2v^{\delta(\sigma + 1)} |cz+d|^{k-2\sigma +K_5}.
\end{eqnarray*}
Choosing $\sigma = (k +K_5)/2$ leads to
\begin{eqnarray*}
g_l(\tau) \leq K_2v^{\delta((k+K_5)/2 + 1)}.
\end{eqnarray*}

\medskip
In the cuspidal case we have $\delta = 0$, whence $g_l(\tau)$ is \emph{bounded} in $\frak{H}$.
Then 
\begin{eqnarray*}
|f_l(\tau)| = y^{-\sigma}g_l(\tau) = O(y^{-k-K_5)/2}).
\end{eqnarray*}
By a standard argument this implies that the Fourier coefficients of $g_l(\tau)$
satisfy $a(n) = O(n^{(k+K_5)/2})$ for $n \rightarrow \infty$. In the holomorphic case there is a similar argument (\cite{KM1}) wherein the exponent is doubled. We have proved

\begin{thm} Suppose that $F(\tau) \in \mathcal{H}(k, \rho)$. There is a constant $\alpha$ depending only on $\rho$ such that the Fourier coefficients of $F(\tau)$
satisfy $a(n) = O(n^{k+ \alpha})$ for $n \rightarrow \infty$. If $F(\tau)$ is cuspidal then
$a(n) = O(n^{k/2+\alpha/2})$ for $n \rightarrow \infty$.
$\hfill \Box$
\end{thm}

\bigskip
\noindent
Marvin Knopp, Department of Mathematics, Temple University, Philadelphia, Pa. 19122\\
Geoffrey Mason, Department of Mathematics, UC Santa Cruz, Ca 95064, USA; gem@cats.ucsc.edu

\end{document}